\documentclass[12pt,letterpaper]{article}
\usepackage[hmargin=1in,vmargin=1in]{geometry}

\usepackage{amsmath}
\usepackage{amsfonts}
\usepackage{amsthm}
\usepackage{amssymb}
\usepackage{array}
\usepackage{caption}
\usepackage{color}
\usepackage{enumitem}
\usepackage{epic,eepic}
\usepackage{graphics}
\usepackage{graphicx}
\usepackage{hyperref}
\usepackage{mathtools}
\usepackage{rotating}
\usepackage{verbatim} 

\newtheorem{prop}{Proposition}

\newtheorem{theorem}[prop]{Theorem}

\newtheorem{conjecture}[prop]{Conjecture}

\begin{document}

\title{Dragons and Kasha}
\author{Tanya Khovanova\\MIT}
\maketitle

\begin{abstract}
Kasha-eating dragons introduce advanced mathematics. The goal of this paper is twofold: to entertain people who know advanced mathematics and inspire people who don't. 
\end{abstract}

\begin{quote}
Suppose a four-armed dragon is sitting on every face of a cube. Each dragon has a bowl of kasha in front of him. Dragons are very greedy, so instead of eating their own kasha, they try to steal kasha from their neighbors. Every minute every dragon extends four arms to the four neighboring faces on the cube and tries to get the kasha from the bowls there. As four arms are fighting for every bowl of kasha, each arm manages to steal one-fourth of what is in the bowl. Thus each dragon steals one-fourth of the kasha of each of his neighbors, while at the same time all of his own kasha is stolen. Given the initial amounts of kasha in every bowl, what is the asymptotic behavior of the amounts of kasha?
\end{quote}

Why do these dragons eat kasha? Kasha (buckwheat porridge) is very healthy. But for mathematicians, kasha represents a continuous entity. You can view the amount of kasha in a bowl as a real number. Another common food that works for this purpose is soup, but liquid soup is difficult to steal with your bare hands. We do not want to see soup spilled all over our cube, do we? If kasha seems too exotic, you can imagine less exotic and less healthy mashed potatoes.

How does this relate to advanced mathematics? For starters, it relates to linear algebra \cite{S}. We can consider the amounts of kasha as six real numbers, as there are six bowls, one on each of the six faces of the cube. We can view this six-tuple that represents kasha at each moment as a vector in a six-dimensional vector space of possible amounts of kasha. To be able to view the amounts of kasha as a vector, we need to make a leap of faith and assume that negative amounts of kasha are possible. I just hope that if my readers have enough imagination to envision six four-armed dragons on the faces of the cube, then they can also imagine negative kasha. The bowl with $-2$ pounds of kasha means that if you put two pounds of kasha into this bowl, it becomes empty. For those who wonder why dragons would fight for negative kasha, this is how mathematics works. We make unrealistic assumptions, solve the problem, and then hope that the solution translates to reality anyway.

Back to the dragons. After all the kasha is redistributed as a consequence of many arms fighting and stealing, the result is a linear operator acting on our vector space, which we will call the \textit{stealing operator}. The 6-by-6 matrix of the stealing operator depends on how we number the faces of the cube. For the numbering in my head it looks like this:
\[
A= 
\begin{bmatrix}
0 & 1/4 & 1/4 & 1/4 & 1/4 & 0 \\
1/4 & 0 & 1/4 & 0 & 1/4 & 1/4 \\
1/4 & 1/4 & 0 & 1/4 & 0 & 1/4 \\
1/4 & 0 & 1/4 & 0 & 1/4 & 1/4 \\
1/4 & 1/4 & 0 & 1/4 & 0 & 1/4 \\
0 & 1/4 & 1/4 & 1/4 & 1/4 & 0 \\
\end{bmatrix}
\]

No matter how you number the faces, all numbers in the matrix are equal to either 0 or $1/4$ because dragons take 0 from themselves and from the opposite dragons and one fourth of the kasha from the other dragons. Exactly four numbers in each row and column should be $1/4$, because this is the number of neighbors the kasha is stolen from, as well as the number of neighbors the kasha from one bowl goes to. The diagonal must have all zeros, because the dragons do not steal from themselves.

To calculate how kasha redistributes after the first fight, we can take the initial distribution of kasha as our vector and multiply it by our matrix. To see what happens after many steps we need to multiply by matrix $A$ many times. Or, we can find powers of matrix $A$ first and then apply this to our vector. The beauty of finding the powers first is that we can see how the stealing operator transforms after many applications without knowing the initial distribution. 

As a true mathematician, I am lazy. I do not want to multiply matrices. I do not even want to type them into my calculator. I want to solve this problem by using my knowledge and the power of my brain without getting off my couch.

But how do I find the asymptotic distribution without multiplying the matrices many times? Mathematicians have found a way to quickly calculate powers. The idea is a \textit{diagonalization} of a matrix. Suppose I find an invertible matrix $S$ and a diagonal matrix $\Lambda$ such that $A=S\Lambda S^{-1}$. Then powers of $A$ are: $A^n= S \Lambda^n S^{n-1}$. 

I am getting excited. It is really easy to compute powers of a diagonal matrix. If $\Lambda$ has $\{\lambda_1,\lambda_2,\ldots,\lambda_6\}$ on the diagonal, then its $n$-th power is again a diagonal matrix with $\{\lambda_1^n,\lambda_2^n,\ldots,\lambda_6^n\}$ on the diagonal. If some of the lambdas are less than 1 in absolute value, their powers tend to zero when $n$ increases, and this is what this problem is about: To find and get rid of negligible behavior for large $n$.

So how do I find the lambdas? The fact that matrix $\Lambda$ is diagonal means that for every $1 \leq i \leq 6$ there exists a vector $v_i$ such that $\Lambda v_i = \lambda_i v_i$. Such vectors are easy to find: for example $v_3=(0,0,1,0,0,0)$. We can choose $v_i$ as a vector with all zeros except 1 in the $i$-th place. Now let's get back to $A$: $\Lambda v_i = \lambda_i v_i$ means $A(Sv_i) = \lambda_i (Sv_i)$. We found vectors such that the stealing operator multiplies them by a constant. Such vectors are called \textit{eigenvectors} and the corresponding constants are called \textit{eigenvalues}. We actually didn't ``find'' such vectors yet. We discovered that if matrix $A$ is diagonalizable, such vectors should exist. We can find the diagonalization by finding eigenvalues and eigenvectors. 

Let's find them. Wait a moment! I'm having an attack of laziness again. I do not want to lift a single finger; instead I'd rather use my brain to figure out what these eigenvalues might be. Can an eigenvalue of the stealing operator have an absolute value more than 1? Suppose there is such a value: $Av = \lambda v$, where $|\lambda| > 1$. Consider the luckiest dragon with the largest absolute value of his kasha. As he gets the average of what his neighbors own, his kasha's absolute value can't increase after the fight. On the other hand, the absolute value of his kasha gets multiplied by $|\lambda| > 1$, which is a contradiction. What was this dragon thinking? His energy would have been better spent protecting his kasha rather than stealing.

Let me summarize: to find the limiting behavior we need to find eigenvalues and their corresponding eigenvectors. Eigenvalues with the absolute value more than 1 do not exist. Eigenvalues with the absolute value less than 1 might exist, but in the limit the corresponding vectors are multiplied by zero, which means we might not need to calculate them. Now we need to find eigenvalues with absolute value 1. Let's start searching for an eigenvector with an eigenvalue that is exactly 1. Now that I think about it, if every dragon has one pound of kasha, their fighting is a complete waste of time: it doesn't change anything. After the fight, each dragon will have one pound of kasha. In other words, vector $(1,1,1,1,1,1)$ is an eigenvector with eigenvalue 1. Such kasha distribution that doesn't change from fight to fight is called a \textit{steady state}. If all other eigenvectors have absolute values less than 1, then $\Lambda^\infty$ has exactly one 1 on the diagonal. As such, it is a matrix of rank 1, which means $A^\infty$ is also rank 1 matrix. Therefore, the asymptotic behavior is proportional to one particular distribution which has to be the steady state. So far, I haven't done any calculations, but I do have a conjecture:  

\begin{conjecture}
Asymptotically every dragon gets $1/6$ of the total amount of kasha. 
\end{conjecture}

The value of $1/6$ comes from the fact that the total amount of kasha doesn't change during stealing process.

Now let's see if I can prove my conjecture. Matrix $A$ looks quite special. Maybe there is something about it we can use. One might recognize this matrix as a Markov matrix of a random walk. To elucidate, let me define these new words. A \textit{Markov matrix} is a matrix with non-negative elements such that each column sums to 1. Such matrices describe transitions between states. The matrix elements are probabilities and column $i$ represents the probabilities of transitioning from state $i$ to all the other states. Matrix $A$ is Markov, and it is even more special than that. All non-zero numbers in every column are the same and equal to 1/4. That means in our process a new state is chosen randomly from a list of four states. This process is called a \textit{random walk}.

Where are random walks coming from in our puzzle? The dragons are not moving! Technically their arms are moving, but they are not walking. This might sound crazy, but in this puzzle the kasha is walking. Imagine a tiny piece of kasha. After each fight it moves from one face of the cube to the neighboring face. We can assume that each tiny piece of kasha has a will of its own. A piece of kasha flips a coin, or more precisely, it flips a coin twice, which is equivalent to flipping two coins once. Using the flips this tiny piece of kasha chooses randomly one of the four hands which grabs it. This approach doesn't change our problem. Each dragon still gets one quarter of the kasha from each bowl they are fighting over. 

After each fight each tiny piece of kasha chooses a new bowl randomly. This is its random walk routine: ``walking'' from a bowl to a bowl. The dragon fight is meaningless. In our new setting, the dragons do not control a particular part of kasha they get. The power and the decision making is transferred to kasha.

We want to calculate the probabilities of where each tiny piece of kasha can end up after many steps. This kasha hike seems like a very different problem from our dragon brawl, but the mathematical description is the same. Let me represent the starting position of a tiny piece of kasha as a vector in the six-dimensional space of faces, with 1 marking the face the piece is starting at. To find the probability of where it can be after the first step, we need to multiply the starting vector by matrix $A$: this is the same $A$ that we had for the kasha-fighting dragons. To find the probability distribution of where the piece of kasha can end up after many steps, we need to find the asymptotic behavior of matrix $A^\infty$. How nice! We can solve the dragon-fighting problem and the kasha-walking problem with the same matrix.

Our conjecture translated into the kasha-walking problem states that after many steps the probability of each piece of kasha ending up on a particular face is $1/6$. It is uniform and doesn't depend on the starting face. So, what does the theory of Markov processes and random walks says about my conjecture? The theory \cite{B} says the following:

\begin{theorem}
The steady state is the limiting distribution if the process is irreducible and aperiodic.
\end{theorem}

Wait a minute. Allow me to explain the two new words in my theorem. \textit{Irreducible} means that the tiny piece of kasha can reach any face of the cube. Our process is irreducible because all the faces are connected. An irreducible Markov chain is \textit{aperiodic} if the piece of kasha is able to walk to a particular face at irregular time intervals without periodicity restrictions. One of the ways to prove the aperiodicity of our process is to show that the kasha piece can, after more than one step, end up at any face of the cube of its choosing.  As I'm still lounging on my couch, I'll leave the proof up to you. 

Anyway, we see that the kasha's random walk is irreducible and aperiodic and therefore tends to its steady state. If the walking kasha ends up on any face with the same probability, then the kasha-fighting dragons will end up in the steady state with the same amounts of kasha.

The conjecture is proven with the reference to an advanced theory and a famous theorem, but I would like to prove it in such a way that the reader can actually check that indeed the steady state has to be the limiting behavior. For this I invoke representation theory. Let us abandon Markov and find a group: representation theory wants a group to represent. We will use the group of rigid motions of the cube. The group acts on the cube, and by extension on the six-tuples of the amounts of kasha. 

This action is called a \textit{representation} of the group. An element $g$ of the group moves the cube with respect to itself. That means it shuffles the six faces of the cube in some way. In this 6-dimensional representation, we assign a matrix $A_g$ to the element $g$. The matrix shuffles the amounts of kasha to match the way faces were shuffled by $g$.

Our dragons respect the group action. Each dragon on each face does exactly the same thing. In other words, the stealing operator commutes with any motion of the cube: you can swap stealing kasha with rotating the cube. If dragons steal kasha first and then the cube is rotated, the result is the same as it would be if they had done these actions in the opposite order. An operator that commutes with the action of the group on our vector space is called an \textit{intertwining operator} of this representation. That means our stealing operator is actually an intertwining operator.

Now we are well into representation theory. We have a 6-dimensional representation of our group. This is a lot of dimensions. Can we simplify this representation? The building blocks of any representation are \textit{irreducible representations}. These are the representations that have no nontrivial invariant subspaces. Let's look at the steady state. This is a 1-dimensional invariant subspace. Indeed, if dragons have the same amounts of kasha, after a cube motion the faces will change, but they will still have the same amounts of kasha. We found one building block. The beauty is that our representation decomposes into irreducibles. That is, there is a complementary representation to the steady state. The complementary 5-dimensional invariant subspace is the subspace of kasha such that the total amount of kasha is zero. Clearly, this 5-dimensional representation is invariant: if we move the cube, the total amount of kasha will not change and will stay zero.
 
Fortunately or unfortunately, our 5-dimensional representation is not irreducible. Why do we want irreducible representations anyway? The idea is that they are the smallest building blocks of any representation. That means they are the simplest we can get, and we hope that everything including the intertwining operator will simplify for each of the irreducible representations. For example, our stealing operator is really simple when acting on the steady state: the operator doesn't change the state. The following statement \cite{E} is the reason to try to find irreducible representations:

\begin{theorem}\label{thm:repr}
If a complex representation of a group can be decomposed into non-isomorphic irreducible representations, then the intertwining operator acts as a scalar on each irreducible representation of the group.
\end{theorem}

A complex representation? If you can imagine negative kasha, you ought to be able to imagine an imaginary kasha. So we just assume that the amounts of kasha are complex numbers. That makes our six-tuples a 6-dimensional complex vector space and our representation a complex representation.
Now we need to continue decomposing.  

The cube has a natural mirror symmetry that swaps the amounts of kasha on the opposite faces of the cube. That means we can decompose the 6-dimensional space of amounts of kasha into two 3-dimensional subspaces that do not change after any rotation: the first subspace has the same amounts of kasha on the opposite faces, and the second subspace has the opposite amounts of kasha on the opposite faces. The 3-dimensional subspace that has the same amounts of kasha on the opposite faces contains a 1-dimensional irreducible representation with the same amounts of kasha on every face. That means we can decompose this 3-d representation into two: a 1-dimensional one we already know about and its complement.

So far we have decomposed the 6-dimensional vector space into the following three representations:

\begin{itemize}
\item 1-dimensional. Every dragon has the same amount of kasha.
\item 2-dimensional. Dragons on the opposite faces have the same amounts of kasha and the total amount of kasha is zero.
\item  3-dimensional. Dragons on the opposite faces have the opposite amounts of kasha.
\end{itemize}

If these representations are irreducible, then they have to be non-isomorphic as they have different dimensions. In this case the stealing operator will act like a multiplication by a scalar. Even if these representations are not irreducible, the stealing operator can still act as a scalar. 

In any case, now it is time to calculate how the stealing operator acts on each representation. 

\begin{itemize}
\item Every dragon has the same amount of kasha. The stealing operator acts as identity.
\item Dragons on the opposite faces have the same amounts of kasha and the total amount of kasha is zero. Consider a red dragon and a blue dragon opposite him. Their four neighbors have the total amount of kasha equal to what the red and the blue dragons have together. That means the neighbors of the red dragon have $-2$ times the amount of kasha the red dragon has. The stealing operator acts as multiplying by $-1/2$.
\item Dragons on the opposite faces have the opposite amounts of kasha. Each dragon is stealing from two pairs of dragons that are opposite each other. The total of the kasha of the neighbors of one dragon is zero. The stealing operator acts as zero. After all this fighting each dragon gets zero kasha. How unproductive.
\end{itemize}

Now we know exactly what happens each time, and we see that asymptotically the stealing operator tends to zero on the two larger invariant subspaces. That means, asymptotically every dragon will have the same amount of kasha. And to tell you a secret, these three representations are indeed irreducible. What I like about this method is that we do not have to believe Theorem~\ref{thm:repr}. We just act on it and get the answer. On top of that, we now know more than the problem asked: how fast we approach the steady state. Hooray to representation theory!

Can we solve the dragon problem without using all these theorems? Yes, we can. For example, here is an elementary solution. By elementary I mean that the most complicated notion it contains is the limit. But if the question asks about the asymptotic behavior, we expect limits anyway. Consider the aftermath of the first fight of our dragons. The dragons on opposite faces get the same amounts of kasha: indeed, they steal equal amounts of kasha from the same dragons who are neighbors to both of them. Now we can assume that dragons on opposite faces have the same amounts of kasha. Consider three dragons sitting on three faces around one corner of the cube. Suppose they have $a$, $b$, and $c$ amounts of kasha. After the fight they get $\frac{b+c}{2}$, $\frac{a+c}{2}$, and $\frac{a+b}{2}$ amounts of kasha. Suppose that numbers $a$, $b$, and $c$ are non-decreasing. That is $a = \min(a,b,c)$ and $c=\max(a,b,c)$. Consider the difference between the maximum and the minimum. Before the fight this difference is $c-a$. After the fight the maximum is $\frac{b+c}{2}$ and the minimum is $\frac{a+b}{2}$. That means, the difference is $(a-c)/2$. The difference between the maximum and the minimum reduces by half after each fight. That means asymptotically this difference tends to zero. Asymptotically, all dragons will have the same amount of kasha.

The solution does not seem too complicated. Why did we discuss advanced mathematics? To be fair, I knew the solution using the representation theory first. So I adapted it to give an elementary explanation. I do not know how easy it is to come up with this explanation without the knowledge of advanced mathematics. In any case, the advanced methods help us move forward and solve more complex problems.

Now let's see what we've learned and solve another dragon-fighting problem:

\begin{quote}
There are $n$ dragons sitting around an $n$-gon-shaped table. Each two-armed dragon is sitting on one side of the table with a bowl of kasha in front of him.  Every minute every dragon extends two arms to the two neighboring polygon's sides and tries to get kasha from the bowls there. As two arms are fighting for every bowl of kasha, each arm manages to steal one-half of what is in the bowl. Thus each dragon steals one-half of the kasha of each of his neighbors, while all of his own kasha is stolen, too. Given the initial amounts of kasha in every bowl, what is the asymptotic behavior of the amounts of kasha?
\end{quote}

Hey, wait a minute! Why do we call them dragons? We could call them greedy people with bad manners. Following the same path as before, we see that there is a steady state with all the kasha the same. So you might expect that the amounts converge to this state. But if it is this easy, why would I offer this puzzle?
Let's use the powerful methods of representation theory we used before. The group we can use here is the rotation group of the $n$-gon. This group is commutative: if we need to perform several rotations, we can do them in any order. Mathematicians call such a group an \textit{abelian} group. Will the commutativity of the group give us an advantage? The representations of abelian groups are especially simple, as you can see in the following statement \cite{E}:

\begin{theorem}
All irreducible complex representations of an abelian group are one-di\-men\-sional.
\end{theorem}

Let's calculate how irreducible representations of the rotation group look like. Every representation is defined by a vector that is multiplied by a scalar if we rotate the table. That is, it is an eigenvector of the rotation operator. Let's pick a designated person with bad manners, and call him Bob. Eigenvectors are defined up to a scaling parameter, so we can assume that Bob has 1 volume of kasha. Suppose Alice, Bob's right neighbor, has $w$ kasha. Suppose we rotate the table and Bob gets Alice's bowl with $w$ kasha. As this is an eigenvector, rotating by one person multiplies the amounts of kasha by a scalar (an eigenvalue) which must be $w$. From here we can calculate that Alice's right neighbor has $w^2$ kasha, and so on. After we rotate $n$ times, where $n$ is the total number of people, we get back to Bob and see that Bob has to have $w^n$ kasha. That means $w^n=1$. Thus, $w$ is a root of unity and for each such root we have an irreducible representation of our group of motions of the table.

How does the stealing operator act on this 1-dimensional representation? By Theorem~\ref{thm:repr} our vector is an eigenvector of the stealing operator. To find the eigenvalue, let's look at Bob and his two neighbors. Bob has $w$ kasha on the right and $w^{-1}$ on the left. So after the first round he will have $(w+w^{-1})/2$ kasha in his own bowl. This is our multiplication coefficient: after a fight every person with bad manners gets his/her kasha multiplied by this number. Given that $w = e^{2\pi ik/n}$, for $0 \leq k < n$, we get that the kasha multiplies by $(e^{2\pi ik/n}+e^{-2\pi ik/n})/2= \cos 2\pi k/n$. 

We are interested in the asymptotic behavior. If the absolute value of the cosine is less than 1, then asymptotically after many iterations we get zero. Suppose the absolute value of the cosine is 1. This can only happen in two cases. For $k=0$, the cosine is 1. In this case $w=1$ and this is our steady state. If $n$ is odd, everything converges to this steady state. If $n$ is even and $k=n/2$ we get another possibility of the absolute value being 1. In this case we have $w=-1$. If Bob's kasha was 1 at the starting point, it will become $-1$ after the first fight. It will continue to fluctuate indefinitely between $1$ and $-1$. Thus we have two eigenvectors for even $n$ that survive the threat of time. The limiting behavior is 2-dimensional.

To summarize, when $n$ is odd, the amounts of kasha converge to the same number for every greedy person with bad manners. If $n$ is even, the amounts of kasha converge to two numbers $a$ and $b$, alternating between people. Asymptotically, after every stealing, the amounts of kasha of one bad-mannered person fluctuate between $a$ and $b$.

Do you remember the discussion on Markov matrices and random walks? As before we can convert this problem to a random kasha-walk on an $n$-gon. What is different here is that when $n$ is even, the process is not aperiodic. A tiny piece of kasha can walk to some of the sides only in an odd number of steps and to other sides in an even number of steps. Thus there is no guarantee of the steady state being the limiting behavior.

After solving these two problems, what can we conclude? That it doesn't pay to be greedy and that mathematics is fun!

\end{document}